\documentclass[letterpaper, 10 pt, conference]{ieeeconf}  

\IEEEoverridecommandlockouts                              

\usepackage{amsmath, amssymb}
\allowdisplaybreaks

\usepackage{amsthm}
\usepackage[dvipsnames]{xcolor}

\usepackage{graphicx} 
\usepackage{subcaption} 
\usepackage{booktabs} 
\usepackage{caption} 
\usepackage{hyperref}
\usepackage{cleveref}
\usepackage{multirow}
\usepackage{makecell}
\usepackage{comment}
\usepackage[ruled,vlined,linesnumbered]{algorithm2e}
\usepackage{cuted}
\SetKwInput{Input}{Input}
\SetKwInput{Output}{Output}
\Crefname{algocf}{Algorithm}{Algorithms}
\Crefname{algocf}{Algorithm}{Algorithms}
\DeclareMathOperator*{\argmin}{arg\,min}

\newcommand{\bfB}{{\bf B}}

\newcommand{\bfZ}{{\bf Z}}

\newcommand{\bfc}{{\bf c}}

\newcommand{\bfx}{{\bf x}}
\newcommand{\bfy}{ {\bf y}}
\newcommand{\bfu}{{\bf u}}

\newcommand{\bfr}{{\bf r}}

\newcommand{\bftheta}{{\boldsymbol \theta}}

\newcommand{\bfgamma}{{\boldsymbol \gamma}}

\newcommand{\bfeta}{{\boldsymbol \eta}}

\newtheorem{theorem}{Theorem}

\title{\LARGE \bf
Zero-Shot Transferable Solution Method for Parametric Optimal Control Problems
}

\author{Xingjian Li, Kelvin Kan, Deepanshu Verma, Krishna Kumar, Stanley Osher and Ján Drgoňa
\thanks{XL and K. Kumar are with Oden Institute, UT Austin; K. Kan and SO are with Department of Mathematics, UCLA; DV is with Department of Mathematical and Statistical Sciences, Clemson University; JD is with Department of Civil and Systems Engineering, Johns Hopkins University.}
\thanks{This work was partially funded by NSF 2208272, 2339678 and 2321040, DOE under award B\&R\# KJ0401010, DE-
SC0S026262, FWP\# CC147,   DOE ASCR SciDAC Institute LEADS, DARPA HR00112590074, and ROSEI at JHU,.}%
}

\begin{document}

\maketitle
\thispagestyle{empty}
\pagestyle{empty}

\begin{abstract}
This paper presents a transferable solution method for optimal control problems with varying objectives using function encoder (FE) policies.
Traditional optimization-based approaches must be re-solved whenever objectives change, resulting in prohibitive computational costs for applications requiring frequent evaluation and adaptation.
The proposed method learns a reusable set of neural basis functions that spans the control policy space, enabling efficient zero-shot adaptation to new tasks through either projection from data or direct mapping from problem specifications.
The key idea is an offline-online decomposition: basis functions are learned once during offline imitation learning, while online adaptation requires only lightweight coefficient estimation.
Numerical experiments across diverse dynamics, dimensions, and cost structures show our method delivers near-optimal performance with minimal overhead when generalizing across tasks, enabling semi-global feedback policies suitable for real-time deployment.

\end{abstract}

\section{INTRODUCTION}
\label{sec:intro}
Optimal control problems arise ubiquitously across engineering disciplines~\cite{bardi1997optimal,cannarsa2004semiconcave,rawlings2017model,MORARI19881}. Although the fundamental mathematical framework remains consistent, practical applications often require solving parametric problems where objectives vary according to task specifications, such as target locations in trajectory planning, terrain characteristics in mobile robotics, or process requirements in manufacturing.
Classical local solution methods~\cite{rao2009survey} are relatively fast but must be solved anew for each instance, whereas global solution methods based on the Hamilton--Jacobi--Bellman equations~\cite{fleming2006controlled} are intractable in high dimensions. Machine learning-based approaches aim to bridge this gap and have achieved considerable success~\cite{ruthotto2020machine,liu2017adaptive}; however, they are typically tied to a fixed objective and lack transferability across tasks.

This paper addresses the challenge of efficiently adapting control policies to new objectives without solving each problem instance from scratch. 
The key insight is to approximate the function space of control policies using a function encoder (FE)~\cite{ingebrand2025functionencodersprincipledapproach}, which learns a reusable set of neural network–parameterized functions. 
Policies for new tasks are then expressed as linear combinations of these basis functions, with task-specific coefficients inferred in a zero-shot manner either from limited trajectory data (LS) or directly from the problem specification (operator). 
As such, the computation cost and data requirement for online adaptation are greatly reduced, enabling efficient transfer without sacrificing accuracy.

The main contributions of the work are:
\begin{itemize}
    \item An imitation learning-based framework for parametric optimal control problems that allows for zero-shot generalization to unseen problem instances without model retraining.  
    \item A semi-global feedback formulation that works for arbitrary inputs and is particularly well-suited when repeated evaluation of the model is required.
    \item Validation through extensive numerical experiments, demonstrating the robustness and near-optimal accuracy on high-dimensional and nonlinear examples.
\end{itemize}

\begin{figure}
    \centering
    \includegraphics[width=\linewidth]{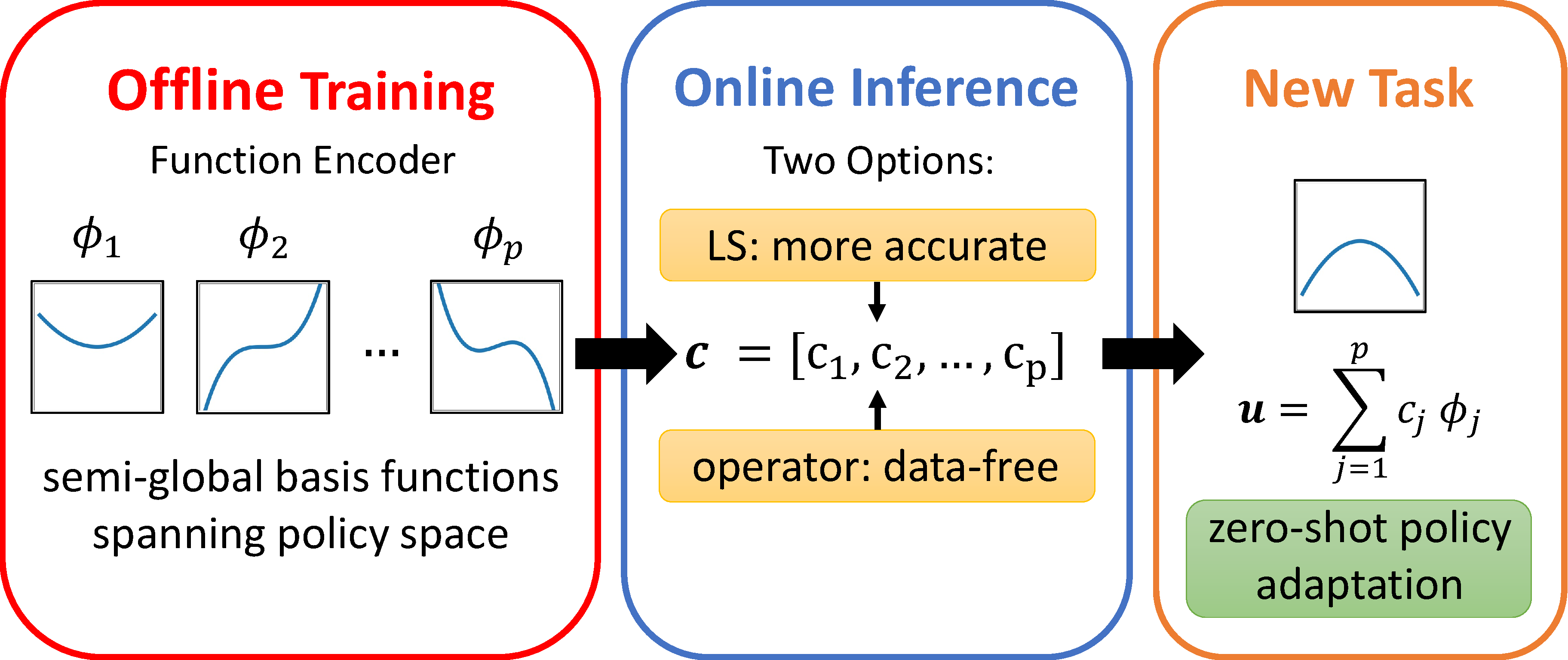}
    \caption{Function  encoder policies: the online–offline decoupling enables efficient and accurate policy adaptation to different optimal control objectives.}
    \label{fig:diagram}
\end{figure}

\section{Background}
\label{sec:background}
This section details the necessary background, covering optimal control and the function encoder framework.

\subsection*{Parametric Optimal Control Problem}
Given an initial state $\bfx(0) = \bfx_0$, we consider a class of optimal control problems where the system dynamics have a fixed form but the objective functional varies with each new task specification. The system evolves according to
\begin{equation}
\dot{\bfx}(t) = f(\bfx(t), \bfu(t), t), \quad \bfx(0) = \bfx_0,
\label{eq:control_dynamics}
\end{equation}
where $\bfx(t) \in \mathbb{R}^n$ is the state, $\bfu(t) \in \mathcal{U} \subseteq \mathbb{R}^m$ is the control input, and $f : \mathbb{R}^n \times \mathbb{R}^m \times [0,T] \to \mathbb{R}^n$ is the dynamics function.
For each task, defined by a conditional variable $\bfeta$, the objective functional in parametric form is
\begin{equation}\label{eq:control_obj_parameterized}
\min_{\bfu \in \mathcal{U}} 
\mathcal{J}(\bfu; \bfeta) 
= \int_0^T L(\bfx(t), \bfu(t), t; \bfeta)\, dt 
+ G(\bfx(T); \bfeta),
\end{equation}
where $L : \mathbb{R}^n \times \mathbb{R}^m \times [0,T] \to \mathbb{R}$ is the running cost, $G : \mathbb{R}^n \to \mathbb{R}$ is the terminal cost, and 
$\mathcal{U}$ is
the space of admissible control functions, which we assume to be sufficiently regular~\cite[Sec. I.3, I.8-9
]{fleming2006controlled}.
For simplicity and without loss of generality, we use $\bfeta$ to denote the conditional variable that captures task-specific variations in $\mathcal{J}$, such as target locations, terrain types, or changes in control penalties. 
In practice, $\bfeta$ can be complex, high-dimensional or implicit. Each choice of $\bfeta$ defines a distinct optimal control problem.

Given the open-loop definition in~\eqref{eq:control_obj_parameterized},
lets consider a closed-loop feedback form~\cite{liberzon2011calculus} given as
\begin{equation}
\bfu(\bfx(t), t; \bfeta) : \mathbb{R}^n \times [0,T] \to \mathcal{U} \subseteq \mathbb{R}^m,
\label{eq:feedback_form}
\end{equation}
for each specific task.
Having the control $\bfu$ to depend on both the current state $\bfx(t)$ and time $t$ allows for control evaluation for any state-time pair, making it more flexible and particularly useful in problems where the initial state $\bfx_0$ can vary. Classical results show that under suitable regularity assumptions, existence of optimal controls in feedback form can be guaranteed; see, e.g., \cite[Sec. I.5]{fleming2006controlled}, or~\cite{Bemporad2002}.

\subsection*{Function Encoder}
Function encoder (FE)~\cite{ingebrand2025functionencodersprincipledapproach} provides a principled framework for representing and transferring tasks in Hilbert spaces by learning a finite set of neural network basis functions. 
Given a Hilbert space $\mathcal{H}$, the FE learns basis functions $\{\phi_1, \phi_2, \dots,\phi_p\}$, parameterized by neural networks, such that any function $h \in \mathcal{H}$ can be approximated as
\begin{equation}
    h(\bfx) \approx \sum_{j=1}^p c_j \phi_j(\bfx; \bftheta_j), 
    \label{eq:fe_formulation}
\end{equation}
for some $\bfc := [c_1, c_2,\dots, c_p]^\top \in \mathbb{R}^p$.
Here we use $\bftheta := \{ \bftheta_1, \bftheta_2, \dots, \bftheta_p \}$ to denote the trainable parameters of the basis functions. This linear representation enables inductive transfer of tasks: new functions are represented by inferring the coefficients $\bfc$ directly from task-specific measurements without retraining the basis functions.

Online computation of coefficients for a given new problem $h$ is obtained by solving a least-squares (LS) problem 
\begin{equation}
    \bfc = \argmin_{\bfc \in \mathbb{R}^p} 
    \| h - \sum_{j=1}^p c_j \phi_j \|_{\mathcal{H}}^2.
\label{eq:ls_problem}
\end{equation}
This admits a closed-form solution given by
\begin{equation}
    \bfc =
    \begin{bmatrix}
        \langle \phi_1, \phi_1 \rangle_{\mathcal{H}} & \cdots & \langle \phi_1, \phi_p \rangle_{\mathcal{H}} \\
        \vdots & \ddots & \vdots \\
        \langle \phi_p, \phi_1 \rangle_{\mathcal{H}} & \cdots & \langle \phi_p, \phi_p \rangle_{\mathcal{H}}
    \end{bmatrix}^{-1}
    \begin{bmatrix}
        \langle h, \phi_1 \rangle_{\mathcal{H}} \\
        \vdots \\
        \langle h, \phi_p \rangle_{\mathcal{H}}
    \end{bmatrix},
    \label{eq:ls_solution}
\end{equation}
where \(\langle \cdot, \cdot \rangle_{\mathcal{H}}\) denotes the inner product associated with the Hilbert space \(\mathcal{H}\), a common choice is the \(L^2\) inner product. 
Both the matrix, known as the Gram matrix, and the right-hand side vector of~\eqref{eq:ls_solution} can be computed directly from task-related measurements or data, using Monte Carlo integration.
In practice, we also incorporate Tikhonov regularization to~\eqref{eq:ls_problem} to improve numerical stability; 
see~\cite{golub1999tikhonov}.
Importantly, FE are supported by a strong theoretical guarantee. 
\begin{theorem}[Universal Function Space Approximation~\cite{ingebrand2025functionencodersprincipledapproach}]\label{thm:UAT}
Let $K \subset \mathbb{R}^n$ be compact and 
let $\mathcal{H} = \{ h : K \to \mathbb{R}^n \mid  \|h\|_{\mathcal{H}} < \infty \}$ be a separable Hilbert space. For any continuous $h \in \mathcal{H}$ and any $\epsilon > 0$, there exist neural network basis functions $\{\phi_1, \phi_2, \dots\}$, some positive integer $P$ and coefficients $\bfc \in \mathbb{R}^P$ such that
\[
  \Big\| h - \sum_{j=1}^{P} c_j\,\phi_j \Big\|_{\mathcal{H}}
  < \epsilon \,\|h\|_{\mathcal{H}}.
\]
\end{theorem}
Here~\Cref{thm:UAT} establishes that, with a sufficient number of basis functions, FEs can approximate any function in $\mathcal{H}$ with arbitrary precision, making them a principled and general-purpose tool for our transfer learning task at hand.

\section{Transferable Solution Method for Parametric Optimal Control Problems}
\label{sec:method}

Many practical optimal control applications require  repeatedly solving problems where system dynamics remain fixed but objectives vary. For example, in trajectory planning, the target destination may change from one instance to another, giving rise to a class of problems that must be solved repeatedly. While  existing approaches (see~\Cref{sec:related_work}) can achieve high accuracy on a fixed problem setting, changes in the objective typically require recomputation of solutions from scratch, leading to substantial computational overhead.

The central challenge hindering model adaptability is to efficiently approximate the family of control policies $\{\bfu^*(\cdot,\cdot;\bfeta)\}_{\bfeta}$ without re-solving each problem instance. To address this, we propose a methodology that directly targets task variability. Our approach leverages the function encoder framework to learn a transferable policy representation, enabling efficient adaptation to new problem specifications with limited, or even no additional data.

\subsection*{Function Encoder Enabled Control Policies}
Our key idea is to approximate the function space of parametric control policies using a function encoder. Specifically, the control policy is modeled as a linear combination of learned basis functions:
\begin{equation}\label{eq:control_fe}
\bfu(\bfx,t;\bfeta) \approx \bfu_{\bftheta}(\bfx,t;\bfeta)  = \sum_{j=1}^p c_j(\bfeta)\,\phi_j(\bfx,t;\bftheta_j),
\end{equation}
where $\{\phi_j(\cdot,\cdot;\bftheta_j)\}_{j=1}^p$ are basis functions realized by neural networks with collective parameters $\bftheta=\{\bftheta_j\}_{j=1}^p$, and $\bfc(\bfeta)$ are coefficients specific to the task $\bfeta$. 

We highlight that the key feature of this formulation is that the basis functions are learned once and are independent of the task parameter $\bfeta$, thereby forming a reusable representation of the policy function space. Consequently, the problem of transferring to a new task reduces to determining the task-specific coefficients $\bfc(\bfeta)$.
This can be achieved in a purely data-driven manner, namely given some dataset $\mathcal{D}$ comprised of $\{\bfx_i, t_i, \bfu(\bfx_i, t_i) \}_{i=1}^{M}$ corresponding to some task $\bfeta$, we consider minimizing the LS formula
\begin{equation}\label{eq:c_min}
    C(\bfc, \mathcal{D}) := \frac{1}{M}\sum_{i=1}^{M} \left\| \bfu(\bfx_i,t_i) - \sum_{j=1}^p c_j\,\phi_j(\bfx_i,t_i) \right\|_2^2.
\end{equation}
over fixed bases $\{ \phi_j\}$ to obtain $\bfc(\bfeta)$. 
In practice, task-dependent datasets ${\mathcal{D}}$ can be obtained through standard open-loop solvers, pretrained policy models, or other means.

\subsection*{Algorithmic Pipeline}
Our methodology follows an offline--online scheme (see~\Cref{fig:diagram}) designed to balance computational efficiency with flexibility. The basis functions are learned only once during the offline phase, while task-specific adaptation takes place online by determining the coefficients $\bfc(\bfeta)$, either from limited observed state-action measurements or via a direct mapping from the task specification~$\bfeta$. This offline-online separation allows the intensive computation to be performed only once during the offline phase, leaving online adaptation lightweight and suitable for real-time control.

\paragraph{Offline Phase} 
FE training is in its core imitation learning. Given task parameterization
$\{\bfeta_1, \dots, \bfeta_N \}$ and their associated datasets
$\{ \mathcal{D}_{S_1},\ldots,\mathcal{D}_{S_N} \}$, we train the basis functions via
Algorithm~\ref{algo:FE_training}. For simplicity, we assume that each dataset contains $M$ labeled observations
$\{((\bfx_i,t_i),\, \bfu_{S_k}(\bfx_i,t_i))\}_{i=1}^M$.

For low-dimensional and structured parameterization of tasks (e.g., target locations), inspired by operator-learning approaches such as \cite{basistobasisoperatorlearning}, we optionally introduce an operator network $\psi:\bfeta \mapsto \bfc(\bfeta)$ with parameters $\bfgamma$. This enables data-free coefficient inference during the online phase. The network is trained via a least-squares reconstruction loss using the fixed bases $\{\phi_j\}$ learned offline; see Algorithm~\ref{algo:operator_training}.

\begin{algorithm}[t]
\caption{FE Training for Control Space}\label{algo:FE_training}
\Input{Learning rate $\alpha$, regularization parameter $\lambda$, task dependent datasets $\{ \mathcal{D}_{S_k}\}_{k=1}^N$, and loss function $C$ from~\eqref{eq:c_min}\;}
Initialize bases $\{\phi_j\}_{j=1}^p$ with params $\bftheta = \{\bftheta_j \}_{j=1}^{p}$\;
\While{not converged}{
    $L \gets 0$\;
  \For {k = 1,2,\dots,N}{
        $\bfc = \argmin\limits_{\bfc} C(\bfc, \mathcal{D}_{S_k})$; \\
    $ L \;\leftarrow\; L \;+\; \| \bfu_{S_k} - \sum_{j=1}^{p} c_j\, \phi_j  \|_{\mathcal{H}}^2$;
  } 
  $
    L \;\leftarrow\; L \;+\; \lambda \sum_{j=1}^{p}  \|\phi_j\|_{\mathcal{H}}^{2}
  $; \\
  $
    \bftheta \;\leftarrow\; \bftheta \;-\; \alpha\, \nabla_{\bftheta}L
  $;
} 
\Output{Trained bases $\{\phi_1, \phi_2\ldots,\phi_p\}$}
\end{algorithm}

\begin{algorithm}[t]
\caption{(Optional) Operator Network Training }
\label{algo:operator_training}
\Input{Trained bases $\{\phi_j\}_{j=1}^{p}$ with parameters $\bftheta$, task--dataset pairs $\{(\bfeta_{S_k}, \mathcal{D}_{S_k})\}_{k=1}^{N}$, learning rate $\beta$, and loss function $C$ from~\eqref{eq:c_min}\; }
Initialize operator net. $\psi:\bfeta\mapsto\mathbb{R}^{p}$ with params $\bfgamma$\;
\While{not converged}{
  $L \gets 0$\;
  \For{k = 1,2,\dots,N}{
  $\bfc = \argmin\limits_{\bfc} C(\bfc, \mathcal{D}_{S_k})$; \\
  $L \;\leftarrow\; L \;+\; \| \bfc - \psi (\bfeta_{S_k}) \|_2^2$;
  }
  $\bfgamma \gets \bfgamma - \beta\,\nabla_{\bfgamma} L$\;
}
\Output{Trained operator network $\psi$}
\end{algorithm}

\paragraph{Online Phase} 
At deployment, the trained basis functions $\{\phi_j\}_{j=1}^{p}$ are fixed, and online adaptation reduces to estimating the task-specific coefficients $\bfc(\bfeta)$ via:
\begin{enumerate}
\item \textbf{Zero-shot LS.} Given some trajectory data for the policy under new task $\bfeta$, we can estimate $\bfc(\bfeta)$ by least-squares projection onto the learned basis functions (minimizing~\eqref{eq:c_min} over measurements).
\item \textbf{Zero-shot operator.} Set $\bfc(\bfeta)=\psi(\bfeta)$ for data-free adaptation using the trained operator network.
\end{enumerate}
We find that in practice, when limited measurements are available for a new task, the LS approach generally gives better performance. On the other hand, online inference is completely data-free using the operator method, though it does incur additional computation cost and requires more training data during the offline phase, and can be particularly hard when $\bfeta$ is high-dimensional or complex. The justification of such trade-off often hinges on each specific problem.

We note that while~\Cref{thm:UAT} guarantees that the bias of the policy adaptation can be made arbitrarily small, 
\eqref{eq:c_min} is a discrete form of~\eqref{eq:ls_problem}, where the inner product in $\mathcal{H}$ is approximated by finite measurements via Monte Carlo integration.
Given some task, denote the minimizers of~\eqref{eq:ls_problem} and~\eqref{eq:c_min} by $\bfc^\star$ and $\bfc$.
They differ due to finite sample approximation of the Gram matrix and right-hand side.
\begin{theorem}\label{thm:LS_Error}
Assume i.i.d.\ sampling of $(\bfx_i, t_i)_{i=1}^{M}$ from the measure defining the inner product in $\mathcal H$. Denote by $\Phi\in\mathbb{R}^{M\times p}$ the matrix whose $i$-th row is $\phi(\bfx_i,t_i)^\top$, $B=\tfrac{1}{M}\Phi^\top\Phi$ is the sample average of the Gram matrix $B^{\star}$ with $\lambda_{\min}(B^\star) > 0$.
Under uniform boundedness assumptions for each $\phi$ and target function $\bfu$, for any $\delta \in (0,1)$, with probability at least $1-\delta$, 
\[
\|\bfc-\bfc^\star\|_2 \le \frac{\lambda_{\min}(B^\star) \sigma_1 + \sigma_2}{\lambda^2_{\min}(B^\star)})\sqrt{\frac{\ln(p/\delta)}{M}}
=\mathcal{O}\big(M^{-1/2}\big),
\]
where $\sigma$ depends on the learned basis functions and $\bfu$.
\end{theorem}
The result follows the concentration bounds for finite sample covariance matrices~\cite{vershynin2018high}, the details of which we provide in the appendix.
Together with~\Cref{thm:UAT}, this implies asymptotic convergence of policy prediction as the number of online samples increase, validating our proposed approach for both offline training and online inference steps.
Empirical validation on function encoder is also available at~\cite{ingebrand2025functionencodersprincipledapproach}.

\section{RELATED WORK}
\label{sec:related_work}
We review approaches to optimal control problems through a local--global lens, then discuss learning-based methods that bridge these extremes, motivating our focus on transferable policy approximation across tasks.

\subsection*{Local Solution Methods (Trajectory Optimization)}
The most common approach for solving optimal control problems is \emph{direct transcription}, which converts the continuous problem into a finite-dimensional optimization problem that is then solved via gradient-based methods~\cite{betts2010practical,rao2009survey}. 
Multiple shooting~\cite{bock1984multiple} extends this idea and often improves numerical stability on long-horizon problems.  
Although highly optimized solvers exist~\cite{Verschueren2021,Andersson2018}, they provide accurate but inherently local solutions that must be recomputed whenever the initial state or objective changes.

\subsection*{Global Solution Methods}
The Hamilton--Jacobi--Bellman (HJB) equation~\cite{fleming2006controlled} characterizes the globally optimal feedback via the value function, from which the optimal control can be recovered as the minimizer of the Hamiltonian. 
This is appealing because it yields a single policy valid across initial conditions given the objective, rather than re-solving a new local trajectory each time. 
In general, solving the HJB equation requires discretizing the state--time domain~\cite{liu1994weighted, osher1991high, achdou2013hamilton}, and the computational and memory costs grow exponentially with the state dimension (the ``curse of dimensionality''), making mesh-based methods non-scalable in practice.
Similarly, solution approaches based on multi-parametric programming~\cite{Bemporad2002,Alessio2009,Kvasnica_2012} face the same issues.
Special cases exist; for instance, LQR problems~\cite{anderson2007optimal, bertsekas2012dynamic} have exact solutions via the Riccati equation; it is powerful yet non-general.

\subsection*{Learning-based Approaches}
ML amortizes solve time into an offline phase to produce feedback policies evaluable online—bridging local trajectory optimization and global HJB solvers. 
Notably, these learned policies are usually \emph{semi-global}: by opting for sampling instead of meshing, they provide reliable feedback over broad regions of the state–time domain, making them scalable, though they remain task-specific in their standard form.
\emph{Data-driven approaches} such as~\cite{hertneck2018learning, karg2020efficient,chen2018approximating, coulson2019data,MADDALENA202011362} train policies on solver-generated trajectories and can achieve fast convergence and high accuracy.
\emph{Model-based and physics-informed approaches}~\cite{han2018solving, onken2022neural, soudbakhsh2023data, drgona2024parametric,Nghiem2023,donti2021enforcing,KORDA2018149,Hewing2020,brunke2022safe,Matschek2023} leverage OC structure (objective terms, HJB residuals, Koopman representation, etc.) to provide robustness and reduce data requirements in training, though the resulting optimization can be considerably more challenging.
\emph{Model-free RL}~\cite{sutton1998reinforcement, hwangbo2019learning} learns feedback without an explicit mathematical model but as a trade-off is often sample-inefficient and may underperform model-based approaches~\cite{verma2025neural}. 
Further integrations of deep learning and optimal control have also yielded promising directions such as~\cite{li2018maximum,onken2021ot,zhang2023mean,kan2025optimal,kan2025stability}, however these discussions are beyond the scope of this paper. 
While ML-based approaches have seen extraordinary success in solving optimal control problems, repurposing a trained model to new tasks efficiently and accurately remains an open challenge. This work aims to address that gap.

\section{Numerical Experiments}
To demonstrate the breadth and generalizability of our approach, we evaluate two categories of optimal control problems differentiated by the source of task variability: changes in the terminal cost and changes in the running cost. The former is illustrated with a linear 2D trajectory planning problem (\Cref{sec:example_2d}) and a 12D quadcopter problem with nonlinear dynamics (\Cref{sec:example_quad}). The latter, motivated by scenarios such as varying terrain and obstacles, is examined through two problems based on a nonlinear bicycle model (\Cref{sec:example_bike}). Together, these examples span low and high state dimensions as well as linear and nonlinear dynamics, providing extensive verification of our method.

Our FE implementation closely follows~\cite{ingebrand2025functionencodersprincipledapproach}. We use multi-head, multi-layer perceptrons (MLPs) for basis function parameterization. Exact training details vary for each example. Default regularization is set to $0.001$.
For data preparation and ground truth comparison we rely on discretization-based open-loop solvers using SciPy or CasADi~\cite{Andersson2018}.
Our code is implemented in PyTorch~\cite{paszke2019pytorch} and can be found \href{https://github.com/XJ-LI/FE_OC.git}{here}. All experiments were conducted on an NVIDIA A100 GPU.

\subsection{2D Path Planning with Different Targets}
\label{sec:example_2d}
We first consider a 2D trajectory planning problem adapted from~\cite{li2024neural} where the goal is to find the optimal path between some initial and target states while avoiding a fixed obstacle. 
The initial state of the agent is not fixed but sampled from a Gaussian distribution, that is, $\bfx_0 \sim \rho = \mathcal{N}((-1.5,-1.5)^\top, 0.4 \cdot \mathbf{I}_2)$. 
We consider linear dynamics where $f = \bfu$ in~\eqref{eq:control_dynamics}, with costs
$$
L(\bfx,\bfu) = \frac{1}{2} \| \bfu \|^2 + Q(\bfx), \  \  G(\bfx) = 50 \cdot \|\bfx - \bfy \|^2,
$$
where $\bfy$ denotes the target state and can change depending on task specification parameter $\bfeta$ in~\eqref{eq:control_obj_parameterized}.
The function $Q(\cdot)$ models an obstacle placed between the initial and target location, given by
$$
Q(\bfx) \;=\; 50 \, \exp(-1.25 \,\|\bfx\|^2 ).
$$
Each task in the example corresponds to a different target state, and our objective is to recover a solution that is evaluable online for any sampled initial state and given target, including those not seen during training.
We prepare the dataset with $16$ different target states over a 2D grid between $1$ and $2$ for both $x$ and $y$. For each target state, we solve the OC problem under $200$ random initializations to obtain trajectories using a direct approach. Here, we set $T = 1$ and use $20$ time steps for the forward integration of the dynamics.

\begin{figure}[!t] 
    \centering
    \begin{subfigure}[b]{0.32\linewidth} 
        \includegraphics[width=\linewidth]{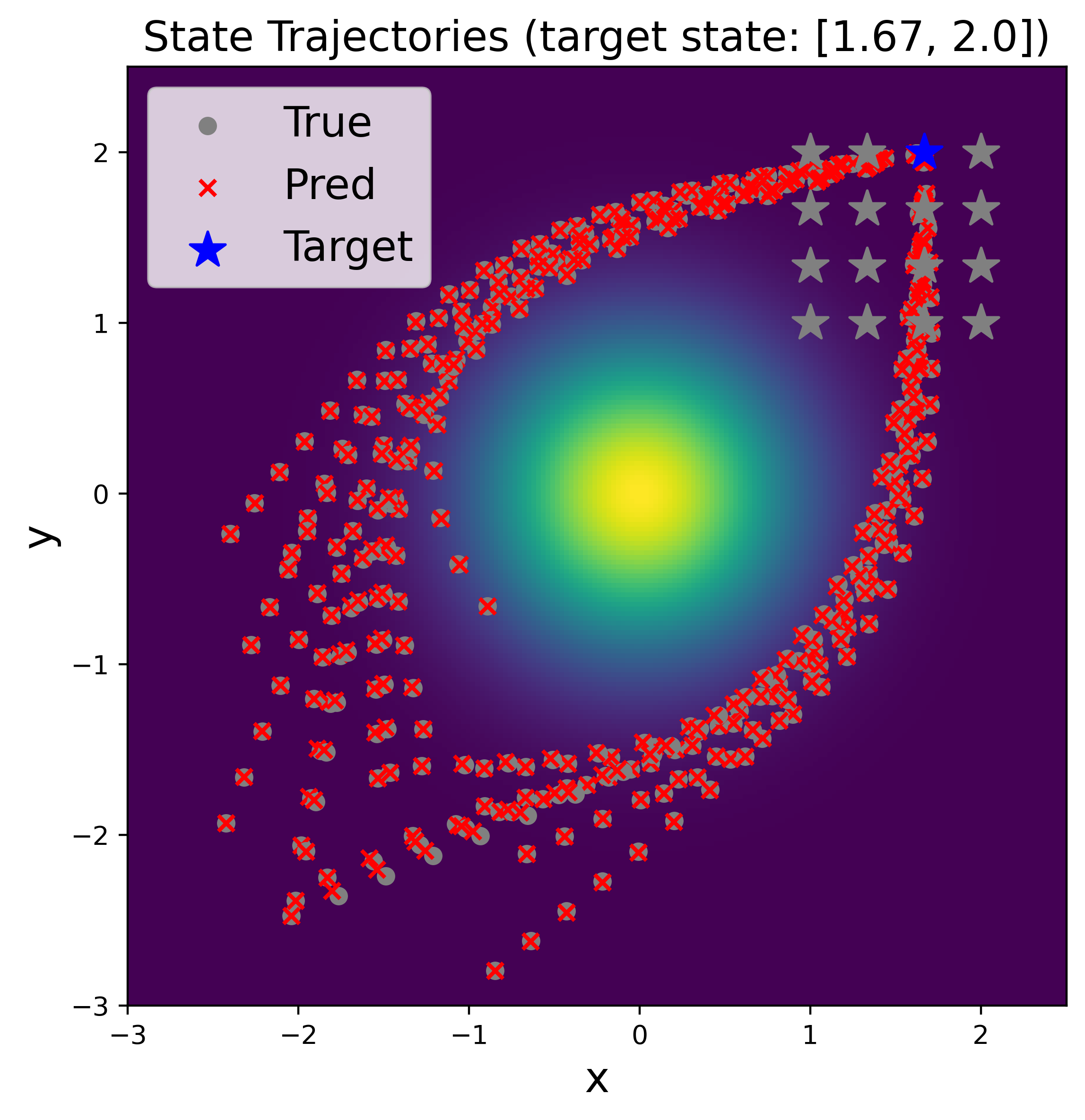}
        \caption{seen target in training, new initial states}
        \label{fig:2d_oc_seen}
    \end{subfigure}
    \hfill 
    \begin{subfigure}[b]{0.32\linewidth}
        \includegraphics[width=\linewidth]{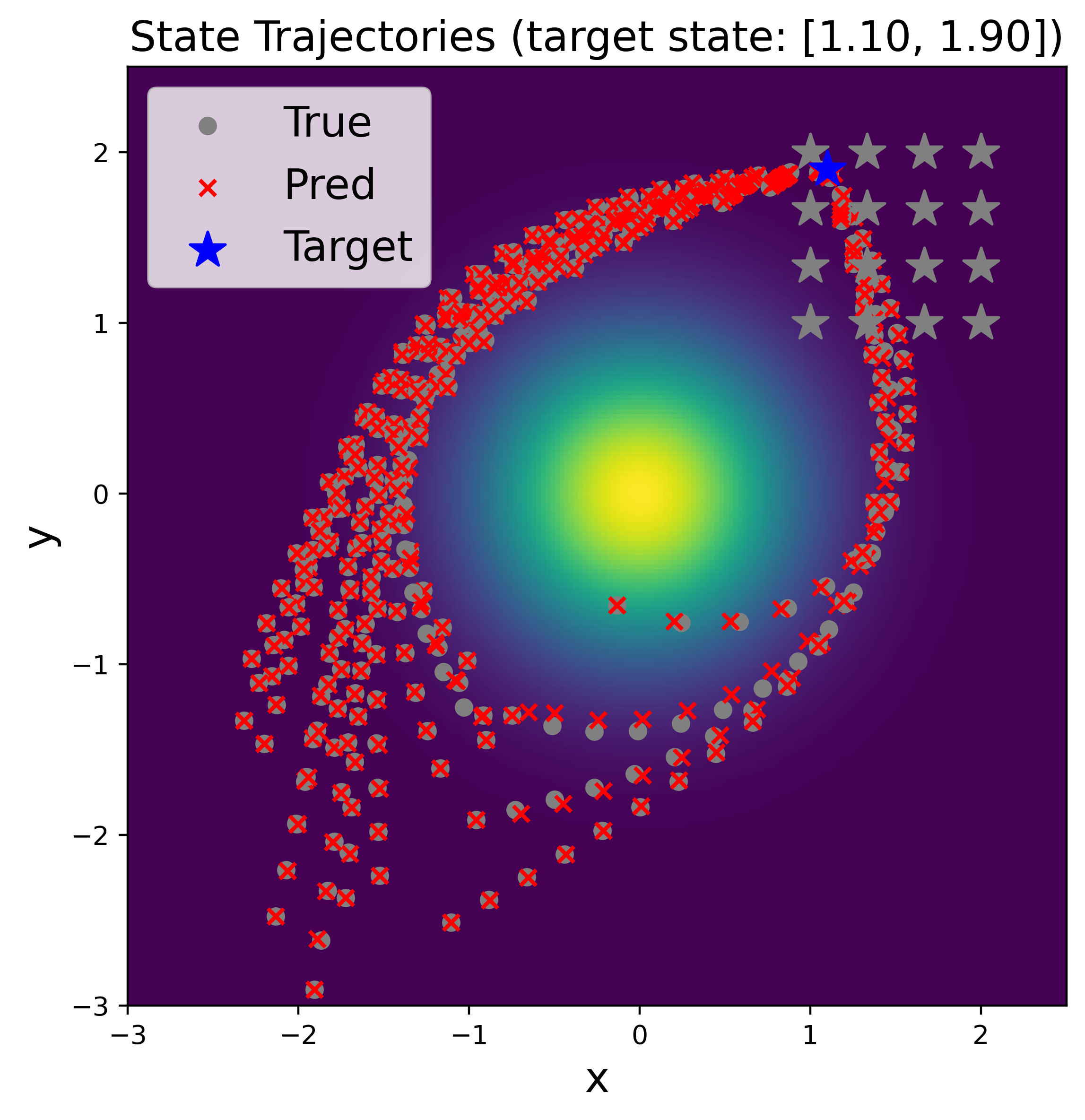}
        \caption{new target, new initial states (interpolation)}
        \label{fig:2d_oc_interp1}
    \end{subfigure}
    \hfill 
    \begin{subfigure}[b]{0.32\linewidth}
        \includegraphics[width=\linewidth]{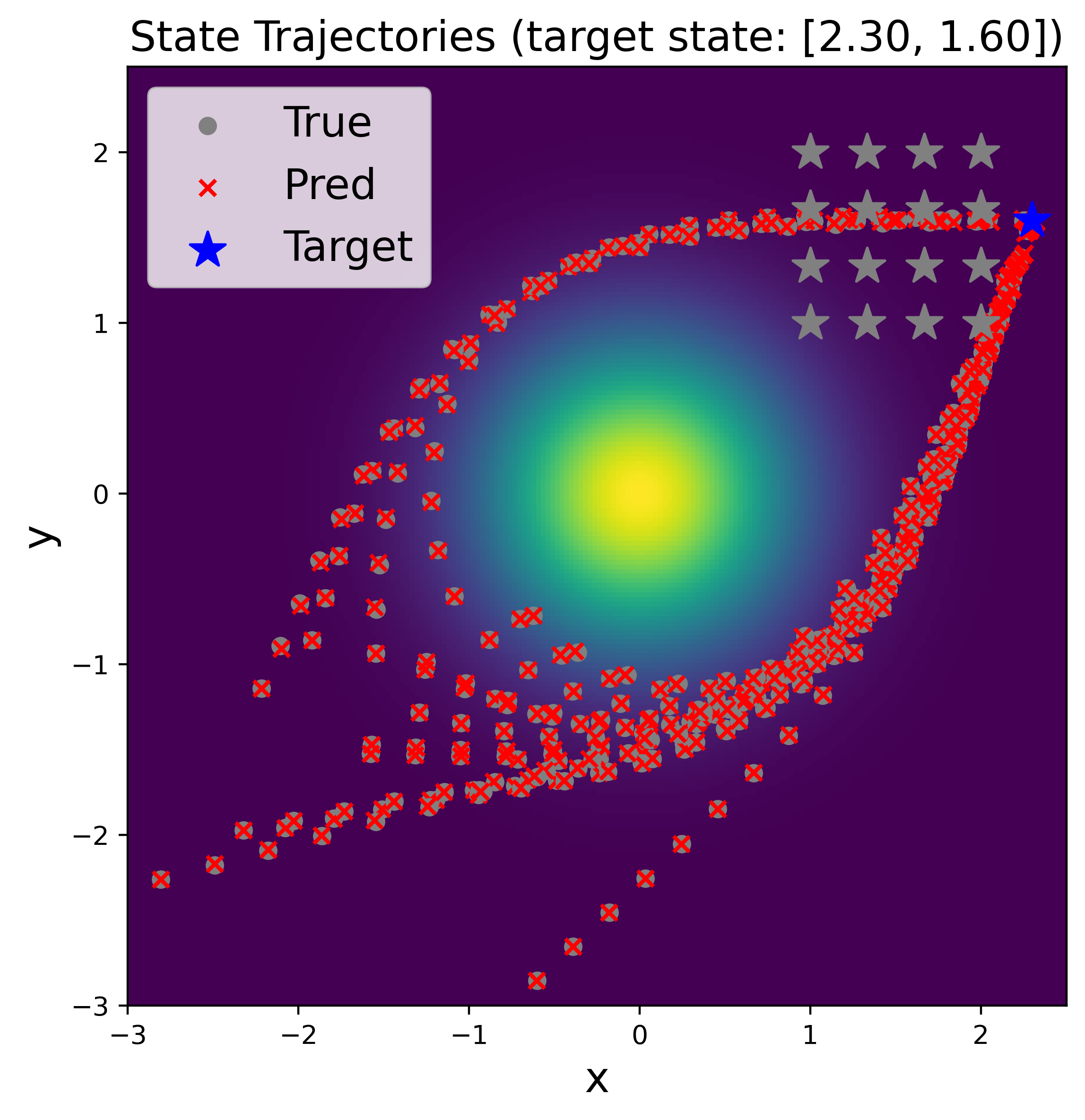}
        \caption{new target, new initial states (extrapolation)}
        \label{fig:2d_oc_extrap}
    \end{subfigure}
    \caption{Generalization results for 2D trajectory planning. All cases test on new initial states, demonstrating semi-global policies that work on both seen and unseen target 
scenarios.}
    \label{fig:2d_oc}
\end{figure}

For the FE parameterization, we use a $4$-layer MLP with hidden size $256$, with $100$ basis functions in total. For training, we use the Adam optimizer~\cite{kingma2014adam} with a learning rate of $0.001$, for a total of $20\text{K}$ steps to ensure convergence.
We present the numerical results in \Cref{fig:2d_oc}, including both the test results on $\bfy$ seen in the training dataset and cases beyond the training domain.
We highlight that in all the tested scenarios, the learned model approximates the ground-truth control accurately without any retraining. To quantify the performance of the learned model, we also report the objective functional value as a reference in \Cref{tab:2d_oc_obj}. Notice that our model retains a low objective functional value with errors below $4\%$ across all test cases.
Note that in the original FE work~\cite{ingebrand2025functionencodersprincipledapproach} model performance is most reliable when generalizing to tasks within the convex hull of the training task distribution, here we also observe promising extrapolation results, such as in~\Cref{fig:2d_oc_extrap}.

Additionally, we present quantitative comparisons between different coefficient estimation methods in \Cref{tab:2d_oc_obj}. The LS approach achieves higher accuracy during online evaluation. Although the operator approach requires additional data and computation for offline training, it offers the advantage of minimal computational cost during the online phase.

\begin{figure*}[!t]  
    \centering
    \includegraphics[width=0.95\textwidth]{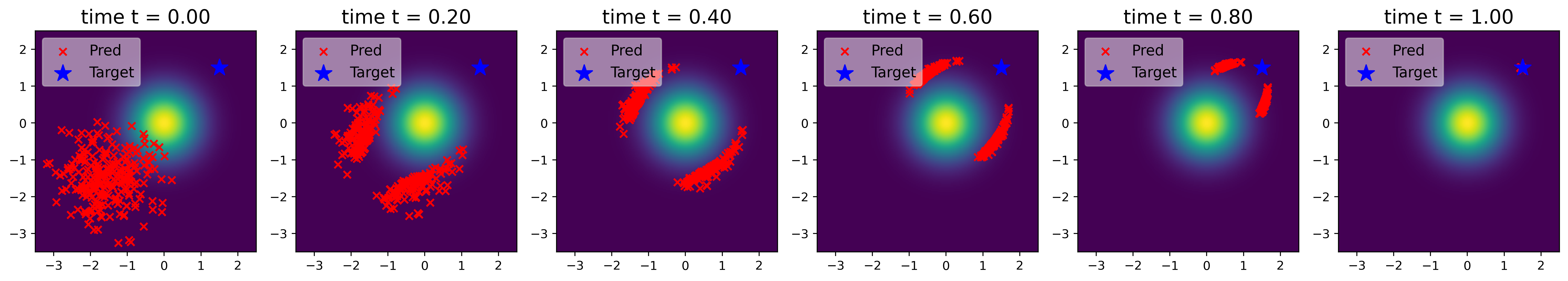}
    \caption{Visualization of the learned control policy. We generate a sample of size $256$ following the specified distribution for the initial state. The plot illustrates how they traverse over time following the learned control, demonstrating consistent performance across the state space.}
    \label{fig:2d_oc_global}
\end{figure*}

We highlight that our approach is fundamentally different from any local solution methods, as the feedback policy can be applied for arbitrary inputs once trained, demonstrated in~\Cref{fig:2d_oc_global}.

\begin{table}[htbp]
\centering
\caption{Average true and predicted objective loss across evaluations for the 2D path planning problem. The average is taken across different target states and initial conditions.}
\begin{tabular}{@{}cccc@{}}
\toprule
Evaluation & \makecell{True \\ Objective Loss} & \makecell{Predicted \\Objective Loss} & \makecell{Inference \\ Method} \\
\midrule
seen target           & 4.8609 & 4.8737 & \multirow{3}{*}{LS} \\
new target interpolation   & 4.8295 & 4.8608 &                       \\
new target extrapolation   & 4.8782 & 5.0751 &                       \\
\midrule
seen target           & 4.8609 & 4.9625 & \multirow{3}{*}{operator} \\
new target interpolation   & 4.8295 & 4.9578 &                       \\
new target extrapolation   & 4.8782 & 5.3214 &                       \\
\bottomrule
\end{tabular}
\label{tab:2d_oc_obj}
\end{table}

\subsection{Quadcopter Path Planning with Different Targets}
\label{sec:example_quad}
We consider controlling a quadcopter under complex dynamics; see~\cite{onken2022neural}. Here, the problem is modeled with a $12$ dimensional state variable $\mathbf{x} \in \mathbb{R}^{12}$ and control inputs consist of a thrust $u$ and torques $(\tau_\psi, \tau_\theta, \tau_\varphi)$.
The objective is to steer the system from a stationary initial state randomly sampled
\begin{equation*}
\begin{aligned}
\mathbf{x}_0 = [\xi_1 -2,\,\xi_2-2,\,\xi_3-2,\,0,\,\dots,\,0]^\top, \ \xi_i \sim \mathcal{N}(0,0.5^2),
\end{aligned}
\end{equation*}
for $i = 1,2,3$, to a target state
$
\mathbf{y} = [y_1, y_2, y_3, 0, \dots, 0]^\top \in \mathbb{R}^{12},
$
where $y_1, y_2, y_3$ can be chosen differently. The system dynamics are given as
\begin{equation}
\begin{cases}
\ddot{x} = \frac{u}{m} \left( \sin(\psi) \sin(\varphi) + \cos(\psi)\sin(\theta) \cos(\varphi) \right) \\
\ddot{y} = \frac{u}{m} \left(-\cos(\psi) \sin(\varphi) + \sin(\psi)\sin(\theta)\cos(\varphi) \right) \\
\ddot{z} = \frac{u}{m} \cos(\theta)\cos(\varphi) - g \\
\ddot{\psi} = \tau_\psi \\
\ddot{\theta} = \tau_{\theta} \\
\ddot{\varphi} = \tau_{\varphi}
\end{cases},
\label{eq:quad_dynamics}
\end{equation}
Here $(x,y,z)$ is the position, $(\psi,\theta,\varphi)$ is the orientation, $g=9.8$ is the gravitational constant, and $m=1$ is the mass. The corresponding first-order system to~\eqref{eq:quad_dynamics} can be easily derived~\cite{onken2022neural}.
We consider the objective function given as
\begin{equation*}
    L(\bfx, \bfu) = \frac{1}{2} \| \bfu \|^2, \  \  G(\bfx) = 500 \cdot \|\bfx - \bfy \|^2,
\end{equation*}
where the penalty on the terminal cost is chosen to be large, following~\cite{onken2022neural}, to ensure that different, even far-away targets can always be reached.

For this example, our goal is to recover transferable solutions with respect to different target states. The high dimensionality of the state space and the strong nonlinearity of the dynamics make the problem particularly challenging.
For data preparation, we use a direct transcription approach for solving trajectory-based solutions given each initial and target states. We fix $T = 2$ and use $50$ time steps for forward integration. For function encoder training, we consider $64$ different target states generated over a uniform grid between $1$ and $4$ in the $x$-, $y$-, and $z$-directions. We then sample initial states and solve for $80$ different trajectories. Additionally, a separate test dataset is prepared with $27$ new tasks; in each case, the test set includes $25$ trajectories from different initial states.

\begin{figure}[htbp]
    \centering
    \begin{subfigure}[b]{0.38\textwidth}
        \centering
        \includegraphics[width=\textwidth]{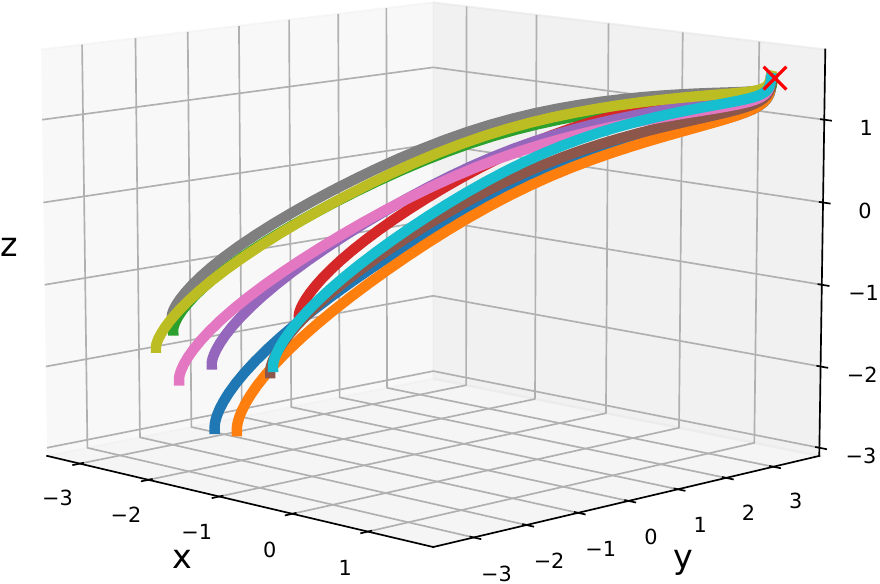}
        \caption{Predictions}
    \end{subfigure}
    \hfill
    \begin{subfigure}[b]{0.38\textwidth}
        \centering
        \includegraphics[width=\textwidth]{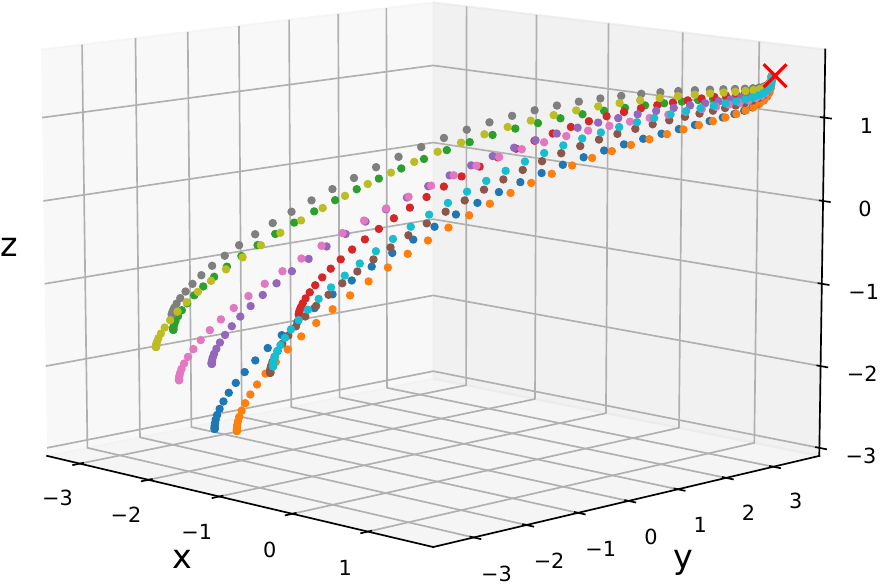}
        \caption{Ground Truth}
    \end{subfigure}
    \caption{Generalization results for the quadcopter path planning problem. Visualization for a new target $\mathbf{y} = [1.5, 3.5, 1.5, 0, \dots, 0]^\top$ not seen during training. We compare model predictions to corresponding true solutions.}
    \label{fig:quad_oc}
\end{figure}

We use the same model architecture and optimization setting as in \Cref{sec:example_2d}, with 100K iterations until convergence.
In~\Cref{fig:quad_oc} we visualize the generalization results of the trained model under a new target state. The model can accurately predict controls and guide the quadcopter to the target location from multiple initializations, demonstrating the effectiveness of our proposed approach.
Quantitative results are displayed in \Cref{tab:quad_oc_obj}. We note that despite the high dimensionality and nonlinearity of the problem, our learned model can achieve high accuracy across different tasks. In fact, under the zero-shot LS inference approach, our policy incurs only $0.4\%$ error in objective value over all $27$ tasks tested. Similar to the findings in \Cref{sec:example_2d}, the LS approach at inference yields more accurate predictions, with additional data requirements as the trade-off.

\begin{table}[htbp]
\centering
\caption{Average true and predicted objective loss across different tasks for the Quadcopter path planning problem, here the average is taken with respect to both different target states and sampled initial conditions. }
\begin{tabular}{@{}cccc@{}} 
\toprule
Evaluation & \makecell{True \\ Objective Loss} & \makecell{Predicted \\Objective Loss} & \makecell{Inference \\ Method} \\
\midrule
seen target     & 276.9546 & 278.8933 & \multirow{2}{*}{LS} \\
new target   & 274.3089 & 275.5074 &                       \\
\midrule
seen target     & 276.9546 & 278.4307 & \multirow{2}{*}{operator} \\
new target   & 274.3089 & 294.5629 &                       \\
\bottomrule
\end{tabular}
\label{tab:quad_oc_obj}
\end{table}

\begin{figure}[h]
    \centering
    \begin{subfigure}[b]{0.3\linewidth} 
        \includegraphics[width=\linewidth]{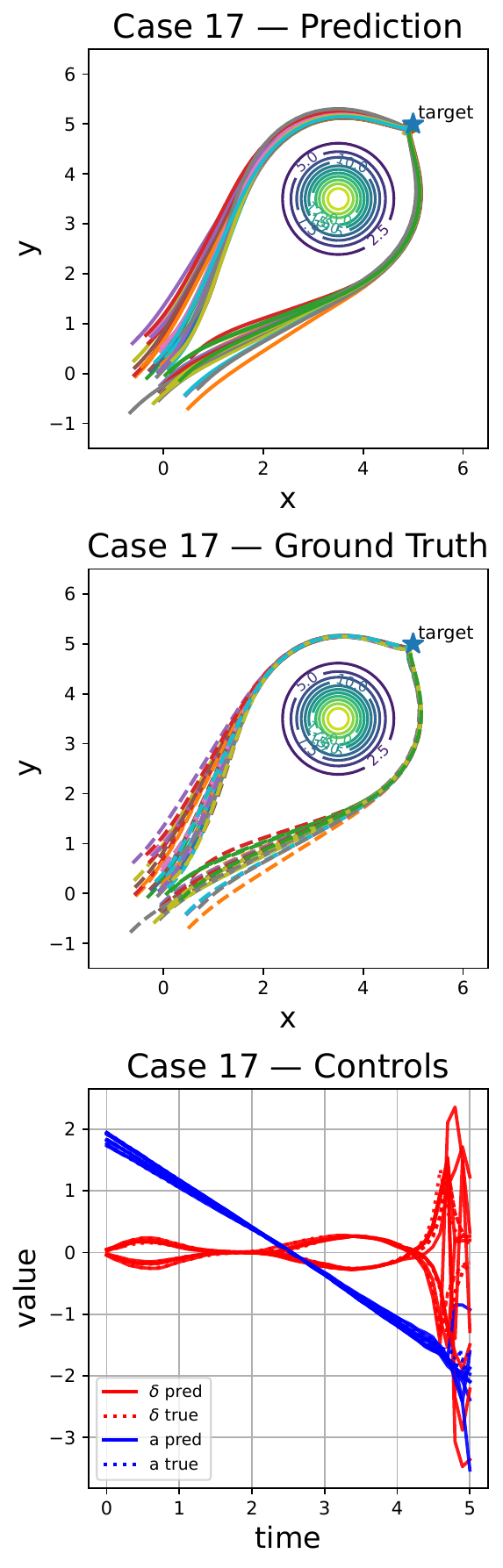}
        \caption{obstacle config 1}
    \end{subfigure}%
    \hfill 
    \begin{subfigure}[b]{0.3\linewidth}
        \includegraphics[width=\linewidth]{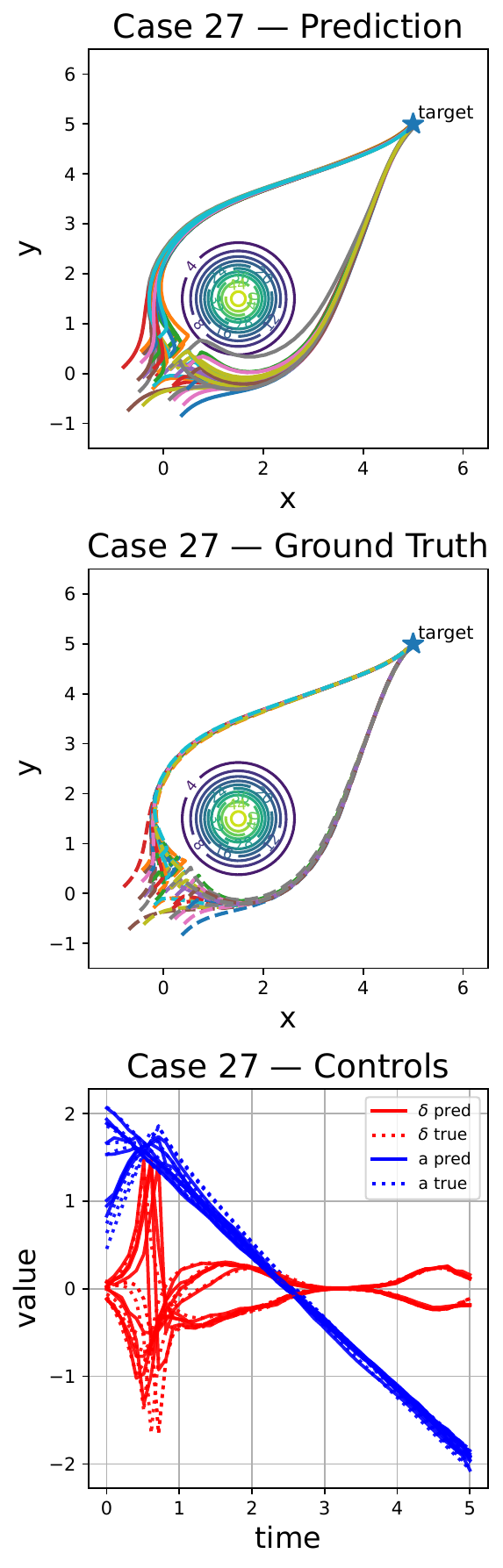}
        \caption{obstacle config 2}
    \end{subfigure}%
    \hfill 
     \begin{subfigure}[b]{0.3\linewidth}
        \includegraphics[width=\linewidth]{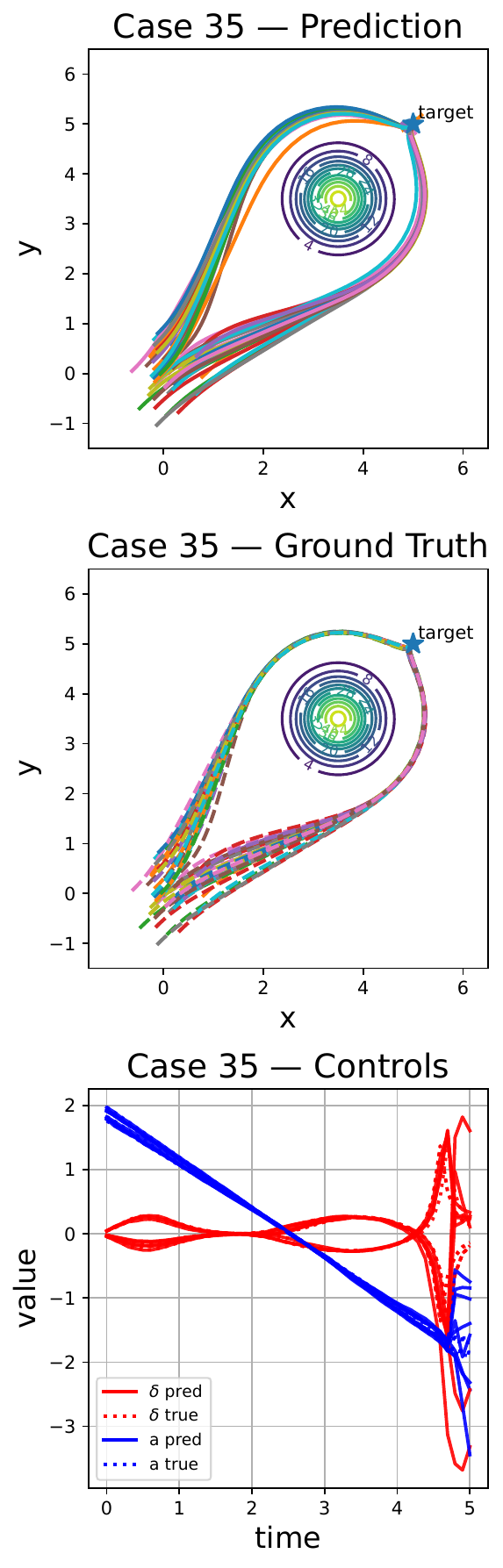}
        \caption{obstacle config 3}
    \end{subfigure}%
    \caption{Visualization of the \emph{three worst performing scenarios} of the learned controller tested over new problem settings for the single obstacle example. Top: Predicted solutions over different initial states. Middle: Ground truth solution calculated w.r.t. the same initial states. Bottom: Visualization of the controls, here for clarity, we only show a few instances.}
    \label{fig:bike_1_gaussian}
\end{figure}

\begin{figure}[h]
    \centering
    \begin{subfigure}[b]{0.3\linewidth} 
        \includegraphics[width=\linewidth]{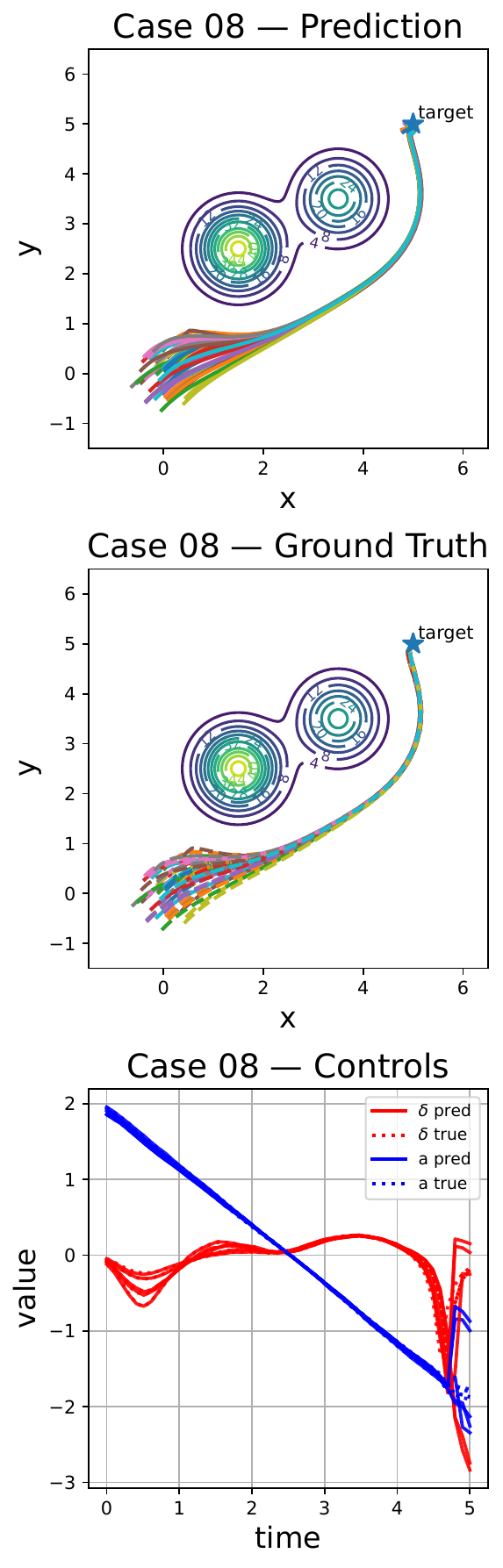}
        \caption{obstacle config 1}
    \end{subfigure}%
    \hfill 
    \begin{subfigure}[b]{0.3\linewidth}
        \includegraphics[width=\linewidth]{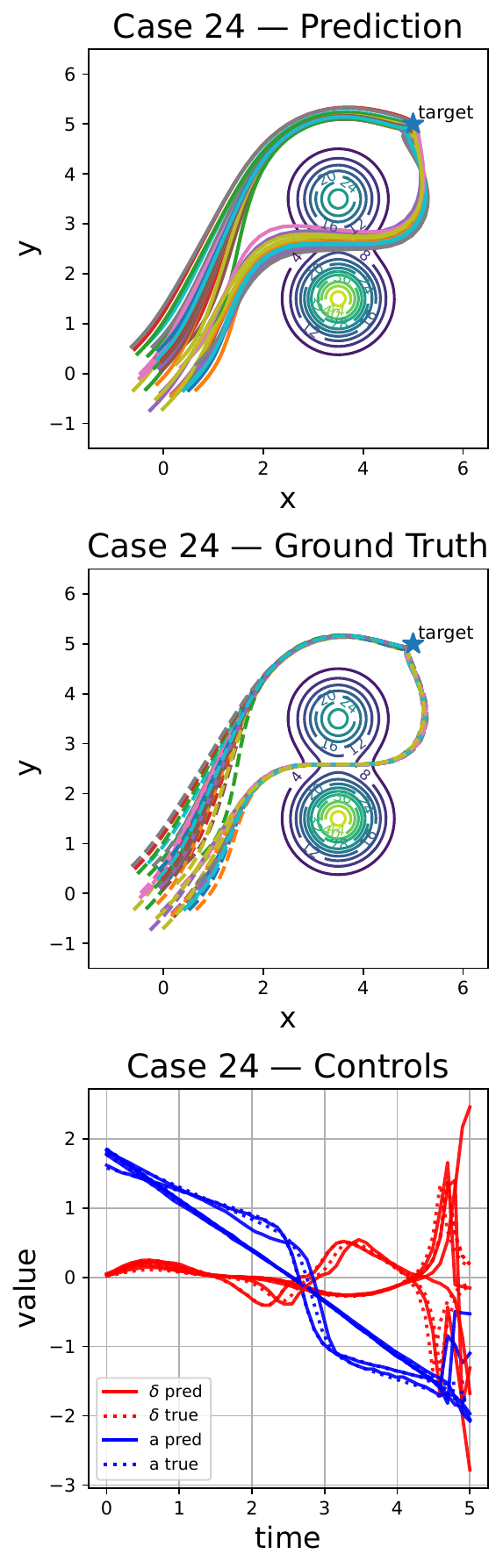}
        \caption{obstacle config 2}
    \end{subfigure}%
    \hfill 
     \begin{subfigure}[b]{0.3\linewidth}
        \includegraphics[width=\linewidth]{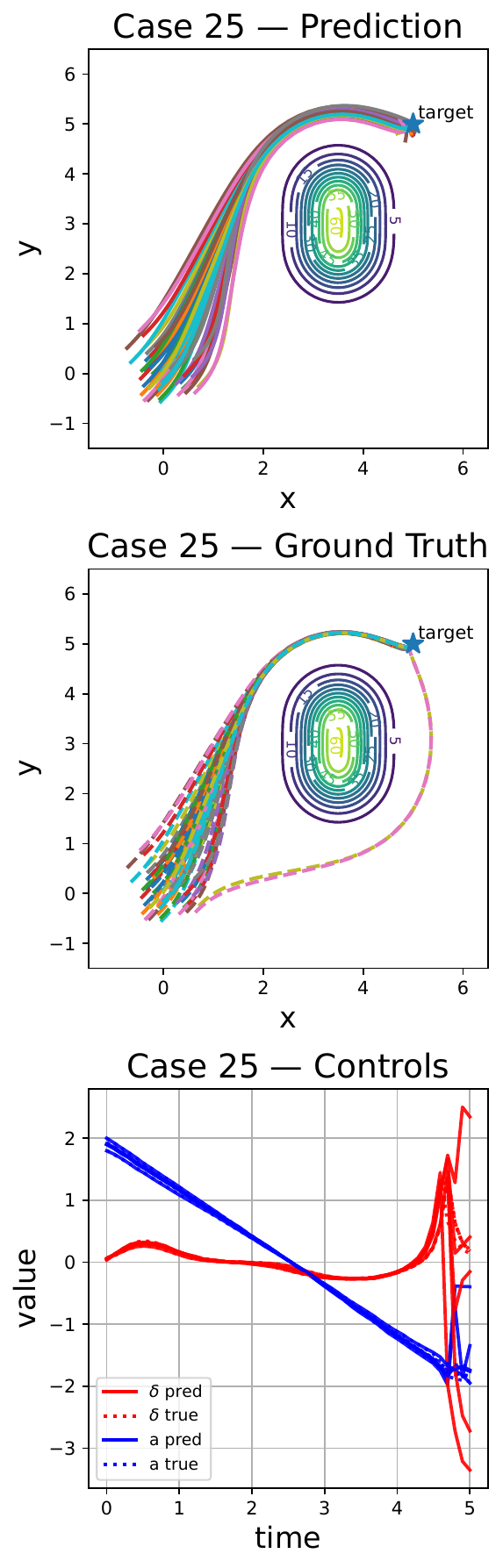}
        \caption{obstacle config 3}
    \end{subfigure}%
    \caption{Visualization of \emph{the three worst performing scenarios} of the learned controller tested over new problem settings for the double obstacle example. Top: Predicted solutions over different initial states. Middle: Ground truth solution calculated w.r.t. the same initial states. Bottom: Visualization of the controls, here for clarity, we only show a few instances.}
    \label{fig:bike_2_gaussian}
\end{figure}

\subsection{Bicycle Control Under Different Obstacle Configurations}
\label{sec:example_bike}
We now consider the scenario where the task specification determines the running cost $L$. This can be common in situations where traversing through different terrain types is needed. Changes to the running cost can substantially alter the control behavior, making these types of problems particularly challenging. Unlike terminal objectives such as target locations, which often admit compact representations, running costs are commonly specified as high-dimensional cost maps over the full domain~\cite{lavalle2006planning}, making them difficult to compress into a form suitable for zero-shot transfer. Hence, in this section, we focus only on the zero-shot LS approach.

For the experiment, we study the motion of a bicycle model. The goal is to find the optimal steering and acceleration to guide the bicycle from a starting point to a target position and orientation over a fixed time horizon. The state of the system at any time $t$ is given by a 4-dimensional vector $\bfx\in \mathbb{R}^4$, written as
$
\bfx(t)=
\begin{bmatrix}
x(t),  y(t), \theta(t), v(t)
\end{bmatrix}^\top,
$
where $(x,y)$ is the 2D position of the bicycle's rear axle, $\theta$ is the heading angle (orientation), and $v$ is the forward velocity.
The control input at $t$ is a 2-dimensional vector $\bfu(t) \in \mathbb{R}^2$ defined as
$
\bfu(t) =
\begin{bmatrix}
\delta(t),   a(t)
\end{bmatrix}^\top,
$
where $\delta(t)$ is the steering angle of the front wheel and $a(t)$ is the acceleration. 
The dynamics of the bicycle are described by the following.
\begin{equation}
    f = 
    \big[ v(t)\cos\theta(t),\; v(t)\sin\theta(t),\; \tfrac{v(t)}{l}\tan\delta(t),\; a(t) \big]^\top,
\label{eq:bike_dynamics}
\end{equation}
where $l=0.5$ is the length of the bicycle's wheelbase.
For this example, while the terminal cost stays quadratic, the running cost consists of penalization on both the control and the state at $t$, which reads
\begin{equation*}
    L(\bfx, \bfu) = \frac{1}{2} \| \bfu \|^2 + Q(\bfx), \  \  G(\bfx) = 50 \cdot \|\bfx - \bfy \|^2,
\end{equation*}
with $Q$ denoting the state penalty given some obstacle.
For this example, we use a time horizon $t \in [0,5]$, the initial state $\bfx_0$ is sampled from a Gaussian distribution centered at the origin, that is,
\begin{equation*}
\bfx_0 = [\xi_1, \xi_2, \frac{\pi}{4}, 0], \quad \xi_1, \xi_2 \sim \mathcal{N}(0,0.35^2).
\end{equation*}
We fix the target state at $\bfy=[5,5,\frac{\pi}{4},0]^\top$ for this example. 

Solving the control problem can be difficult given the highly nonlinear dynamics of the model. For each task, we define a different obstacle configuration that can change the controller behavior, increasing the complexity of the learning problem. 
For demonstration purposes, we model the obstacle using a Gaussian formulation
\begin{align*}
Q(\mathbf{x};A, \boldsymbol{\mu},\sigma)
&= A \,\exp\!\left(-\frac{1}{2\sigma^2}\,\|\mathbf{x}-\boldsymbol{\mu}\|^2\right). 
\end{align*}
Here the parameters $A$, $\boldsymbol{\mu}$ and $\sigma$ control the position and shape of the obstacle and vary from each individual task. 
We generate data for optimal control solutions of $144$ different obstacle settings by randomly selecting $A \in \{30, 40, 50\}$, $\mu_1, \mu_2 \in \{ 1,2, 3, 4 \}$ and $\sigma \in \{0.3, 0.5, 0.7 \}$. For each scenario, we simulate $100$ different trajectories for training.
Another $36$ new cases in obstacle configuration are reserved for testing, each also containing $100$ different trajectories.
We note that for the scenarios where obstacle placement is close to the initial or target states, sharp changes in the ground truth controls are to be expected. We test not only the generalizability but also the robustness of our approach through this example, as the learned model has to account for both smooth and non-smooth policies.

\begin{table}[htbp]
\centering
\caption{Numerical results for the obstacle avoidance example. We display the prediction and ground truth for the cumulative control and obstacle cost, as well as terminal state deviation measured in both $x$ and $y$ direction, all results are evaluated and averaged over settings unseen in training.}
\begin{tabular}{@{}lccc@{}}
\toprule
Evaluation & Control Cost & Obstacle Cost & Term. State Deviation \\
\midrule
\multicolumn{4}{c}{\textbf{Single obstacle}} \\
\midrule
Prediction   & 2.6789 & 0.0463 & 0.0012 \\
Ground Truth & 2.6480 & 0.0385 & 0.00005 \\
\midrule
\multicolumn{4}{c}{\textbf{Double obstacles}} \\
\midrule
Prediction   & 3.0370 & 0.5750 & 0.0046 \\
Ground Truth & 2.9668 & 0.5286 & 0.00004 \\
\bottomrule
\end{tabular}
\label{tab:bike_oc_generalization}
\end{table}

We illustrate the numerical results in \Cref{tab:bike_oc_generalization}. Note that the learned model transfers accurately to new problem settings using limited data for inference. In particular, the controller can learn to avoid the obstacle for arbitrary placement and reach the target state with high precision, demonstrating the capability of our proposed approach.
In \Cref{fig:bike_1_gaussian}, we display the results for $3$ worst performing subproblems out of all tested cases. Notice that even in the worst performing cases, the learned model can provide reasonable and accurate guidance to the system, as demonstrated in the close match between the predicted and ground truth trajectories. 

\subsubsection*{Double  Obstacle Experiment}
To further demonstrate the applicability of our method, we conduct additional experimentation with increasing difficulty. In this example, instead of modeling the terrain using a simple Gaussian obstacle, we introduce a non-overlapping secondary obstacle to the problem. Using the same physics in~\eqref{eq:bike_dynamics},  we define 
\begin{align*}
Q(\mathbf{x}; \{ A_i, \boldsymbol{\mu}_i,  \sigma_i \}_{i=1}^2)
= \sum_{i=1}^2 A_i \exp\!\left(-\frac{1}{2\sigma_i^2}\,\|\mathbf{x}-\boldsymbol{\mu}_i\|^2\right).
\end{align*}
Consider $A_1, A_2 \in \{30, 50\}$, $\sigma_1, \sigma_2 \in \{0.3, 0.5\}$ and the obstacle placement randomly selected over a $3 \times 3$ grid between $1.5$ and $3.5$. We prepare the dataset consisting of $576$ different obstacle configurations, we split them into a $544/32$ training/test split after random shuffling. For each scenario, the dataset contains $100$ trajectories. We note that the secondary obstacle largely increases the difficulty of the problem, as new pathing options lead to different controls; we increase the FE model size to account for this.

We display the corresponding numerical results in \Cref{tab:bike_oc_generalization}, trained model under our approach remains accurate despite the increase in problem complexity.
We display the model predictions and the corresponding true solution in \Cref{fig:bike_2_gaussian}. Here we again only visualize the worst performing cases to show the effectiveness of our approach. We highlight that the new pathing options are correctly captured, and the model performance remains stable even when the underlying ground truth solution exhibits shock-like behavior.

\section{Conclusions}
To summarize, we proposed an FE-based framework for policy transfer in parametric optimal control problems. The framework enables zero-shot adaptation to new tasks with or without data and demonstrates reliable performance across diverse examples. Future work will explore extensions to multi-agent systems with interacting dynamics.

\addtolength{\textheight}{-12cm}   



\bibliographystyle{IEEEtran}
\bibliography{reference}

\makeatletter
\setlength{\textheight}{9.25in} 
\setlength{\topmargin}{-0.25in} 
\makeatother

\onecolumn
\appendices
\section*{APPENDIX}\label{sec:appendix}
\noindent\textbf{Proof of~\Cref{thm:LS_Error}.}\;
We define the population LS coefficients
\[\bfc_* =
    \begin{bmatrix}
        \langle \phi_1, \phi_1 \rangle_{\mathcal{H}} & \cdots & \langle \phi_1, \phi_p \rangle_{\mathcal{H}} \\
        \vdots & \ddots & \vdots \\
        \langle \phi_p, \phi_1 \rangle_{\mathcal{H}} & \cdots & \langle \phi_p, \phi_p \rangle_{\mathcal{H}}
    \end{bmatrix}^{-1}
    \begin{bmatrix}
        \langle \bfu, \phi_1 \rangle_{\mathcal{H}} \\
        \vdots \\
        \langle \bfu, \phi_p \rangle_{\mathcal{H}}
    \end{bmatrix} 
    = \mathbb{E} [\Phi^\top \Phi ]^{-1} \mathbb{E} [\Phi^\top\bfu]
    = \bfB_*^{-1} \bfr_*,
    \label{eq:ls_solution_pop}
\]    
and the corresponding sample average approximation
\[
\bfc =
    \begin{bmatrix}
        \langle \phi_1, \phi_1 \rangle_{M} & \cdots & \langle \phi_1, \phi_p \rangle_{M} \\
        \vdots & \ddots & \vdots \\
        \langle \phi_p, \phi_1 \rangle_{M} & \cdots & \langle \phi_p, \phi_p \rangle_{M}
    \end{bmatrix}^{-1}
    \begin{bmatrix}
        \langle \bfu, \phi_1 \rangle_{M} \\
        \vdots \\
        \langle \bfu, \phi_p \rangle_{M}
    \end{bmatrix} 
    =  [\frac{1}{M}\sum_{i=1}^M \Phi^\top(\bfx_i) \Phi(\bfx_i) ]^{-1} [\frac{1}{M} \sum_{i=1}^M \Phi^\top(\bfx_i)\bfu(\bfx_i) ]
    = \bfB_M^{-1} \bfr_M,
    \label{eq:ls_solution_sa}
\]
where $\langle \cdot, \cdot \rangle_{M}$ is the Monte-Carlo approximation of the inner product $\langle \cdot, \cdot \rangle_{\mathcal{H}}$ using $M$ samples. We consider the following standard assumptions
\begin{enumerate}
    \item $\| \Phi(\bfx) \|_2^2 \le C_1$ for some constant $C_1$ and for all $\bfx$.
    \item $\| \Phi^\top(\bfx) \bfu(\bfx) \|_2 \leq C_2$, for some constant $C_2$ and for all $\bfx$.
\end{enumerate}
Note that we have
\begin{equation}\label{eq:B_star_bound}
    \| \bfB_* \|_2 \leq \| \mathbb{E} (\Phi^\top \Phi)\|_2 \leq \mathbb{E} \left( \| \Phi^\top \Phi\|_2 \right) = \mathbb{E} \left( \|\Phi\|_2^2 \right) \leq \mathbb{E} (C_1) = C_1,
\end{equation}
where we used Jensen's inequality in the second step, and assumption 1 in the second last step. Similarly, we have
\begin{equation}\label{eq:r_star_bound}
    \| \bfr_* \|_2 = \| \mathbb{E} (\Phi^\top \bfu) \|_2 \leq \mathbb{E} \left( \| \Phi^\top \bfu \|_2 \right) \leq \mathbb{E} (C_2) = C_2,
\end{equation}
where we used Jensen's inequality in the second step, and assumption 2 in the second last step.

\noindent Define random matrices
\[
\bfZ_i = \Phi^\top(\bfx_i)\Phi(\bfx_i) - \bfB_*, \quad \text{for} \ i = 1, 2, \dots, M.
\]
Then $\{\bfZ_i\}_{i=1}^M$ are independent, mean-zero random matrices and 
\[
\| \bfZ_i \|_2 \le \| \Phi(\bfx_i)^\top \Phi(\bfx_i) \|_2 + \| \bfB_* \|_2 = \| \Phi(\bfx_i) \|^2_2 + \| \bfB_* \|_2 \le 2C_1,
\]
where in the first step we used triangular inequality, and in the last step we used assumption 1 and~\Cref{eq:B_star_bound}. The variance proxy follows
\begin{align}
\|  \mathbb{E}(\bfZ_i^2) \|_2 &= \| \mathbb{E}(\Phi^\top(\bfx_i)\Phi(\bfx_i)\Phi^\top(\bfx_i)\Phi(\bfx_i) - \Phi^\top(\bfx_i)\Phi(\bfx_i)\bfB_* - \bfB_*\Phi^\top(\bfx_i)\Phi(\bfx_i) + \bfB_*^2) \|_2 \\
&= \| \mathbb{E}(\Phi^\top(\bfx_i)\Phi(\bfx_i)\Phi^\top(\bfx_i)\Phi(\bfx_i)) - \mathbb{E}(\Phi^\top(\bfx_i)\Phi(\bfx_i))\bfB_* - \bfB_* \mathbb{E}(\Phi^\top(\bfx_i)\Phi(\bfx_i)) + \bfB_*^2 \|_2 \\
&= \|  \mathbb{E}(\Phi^\top(\bfx_i)\Phi(\bfx_i)\Phi^\top(\bfx_i)\Phi(\bfx_i)) - \bfB_*^2 \|_2 \\
&\le \| \mathbb{E}(\Phi^\top(\bfx_i)\Phi(\bfx_i))^2 \|_2 + \| \bfB_*^2 \|_2 \\
&\leq \mathbb{E} \| (\Phi^\top(\bfx_i)\Phi(\bfx_i))^2 \|_2 + \| \bfB_*^2 \|_2 \\
\intertext{by Jensen's inequality,}
&= \mathbb{E} \| \Phi^\top(\bfx_i)\Phi(\bfx_i) \|_2^2 + \| \bfB_*^2 \|_2 \\
&= \mathbb{E} \| \Phi(\bfx_i) \|_2^4 + \| \bfB_* \|_2^2 \\
&\leq 2C_1^2,
\end{align}
where in the last step we used assumption 1 and~\Cref{eq:B_star_bound}.
Now applying vector/matrix Bernstein inequalities~\cite{tropp2015introduction}[Theorem~6.6.1], we can get the following bounds with probability $1-\delta$ for $\forall \delta \in (0,1)$, 
\begin{equation}\label{eq:bound_B}
    \|\bfB_M - \bfB_*\|_2 \leq \sqrt{\frac{4 C_1^2  \ln(2p/\delta)}{M}} + \frac{4C_1 \ln(2p/\delta)}{3M},
\end{equation}
where $p$ is the number of basis functions. Following similar steps and we have for the right hand side
\begin{equation}\label{eq:bound_r}
    \|\bfr_M - \bfr_*\|_2 \leq \sqrt{\frac{4C_2^2 \ln((p+1)/\delta)}{M}} + \frac{4C_2 \ln((p+1)/\delta)}{3M}.
\end{equation}

\noindent We then consider the bound for $\| \bfB_M^{-1} - \bfB_*^{-1} \|_2$. Assume sample size $M$ is selected large enough such that 
\begin{equation}\label{eq:N_large_enough} 
    \|\bfB_M - \bfB_*\|_2 \leq \frac{1}{2} \lambda_{\rm min}(\bfB_*).
\end{equation}
Note that $\bfB_M^{-1} - \bfB_*^{-1} = \bfB_M^{-1}(\bfB_*-\bfB_M) \bfB_*^{-1}$, we have
\[
\| \bfB_M^{-1} - \bfB_*^{-1}\|_2 \leq \| \bfB_M^{-1} \|_2 \cdot \| \bfB_*-\bfB_M \|_2  \cdot \| \bfB_*^{-1} \|_2.\label{eq:to_eval}
\]
By Weyl inequality,
\begin{align}\label{eq:weyl}
    \lambda_{\rm min}(\bfB_M) \geq \lambda_{\rm min}(\bfB_*) + \lambda_{\rm min}(\bfB_M - \bfB_*),
\end{align}
also we have
\begin{equation}\label{eq:lmbda_B}
    \lambda_{\rm min}(\bfB_M - \bfB_*) \geq - \| \bfB_M - \bfB_* \|_2.
\end{equation}
By \eqref{eq:weyl} and \eqref{eq:lmbda_B}, we get
\begin{align}\label{eq:bound_BM_inverse}
    \| \bfB_M^{-1} \|_2 &= \frac{1}{\lambda_{\rm min}(\bfB_M)} \leq \frac{1}{\lambda_{\rm min}(\bfB_*) + \lambda_{\rm min}(\bfB_M - \bfB_*)} \leq \frac{1}{\lambda_{\rm min}(\bfB_*) - \| \bfB_M - \bfB_* \|_2} \\
    &\leq \frac{1}{\lambda_{\rm min}(\bfB_*)-\frac{1}{2}\lambda_{\rm min}(\bfB_*)} = \frac{2}{\lambda_{\rm min}(\bfB_*)},
\end{align}
where in the second last step, we used \eqref{eq:N_large_enough}. Continuing evaluating $\| \bfB_M^{-1} - \bfB_*^{-1}\|_2$, and using \eqref{eq:bound_B}, \eqref{eq:bound_BM_inverse}, and that $\| \bfB_*^{-1} \|_2=1/ \lambda_{\rm min}(\bfB_*)$, we have
\begin{align}\label{eq:bound_BM_B_inverse}
    \| \bfB_M^{-1} - \bfB_*^{-1}\|_2 &\leq \| \bfB_M^{-1} \|_2 \cdot \| \bfB_*-\bfB_M \|_2  \cdot \| \bfB_*^{-1} \|_2 \\
    &\leq \frac{2}{\lambda_{\rm min}(\bfB_*)} \left( \sqrt{\frac{4C_1^2 \ln(2p/\delta)}{M}} + \frac{4C_1 \ln(2p/\delta)}{3M} \right) \frac{1}{\lambda_{\rm min}(\bfB_*)} \\
    &= \frac{2}{\lambda_{\rm min}^2(\bfB_*)} \left( \sqrt{\frac{4C_1^2 \ln(2p/\delta)}{M}} + \frac{4C_1 \ln(2p/\delta)}{3M} \right).
\end{align}
In order to bound $\| \bfc - \bfc_* \|_2$, we use the fact
\begin{align}
    \| \bfc - \bfc_* \|_2 &= \| \bfB_M^{-1} \bfr_M-\bfB_*^{-1} \bfr_* \|_2 \\
    &= \| \bfB_M^{-1} \bfr_M - \bfB_M^{-1} \bfr_* + \bfB_M^{-1} \bfr_* -\bfB_*^{-1} \bfr_* \|_2 \\
    &\leq \| \bfB_M^{-1} \bfr_M - \bfB_M^{-1} \bfr_* \|_2 + \| \bfB_M^{-1} \bfr_* -\bfB_*^{-1} \bfr_* \|_2 \\
    &\leq \|\bfB_M^{-1}\|_2 \| \bfr_M - \bfr_* \|_2 + \| \bfB_M^{-1}-\bfB_*^{-1} \|_2 \| \bfr_*\|_2, \\
    & \intertext{here, using \eqref{eq:bound_BM_inverse}, \eqref{eq:bound_r}, \eqref{eq:bound_BM_B_inverse} and \eqref{eq:r_star_bound}, we obtain} \nonumber\\
    & \leq \frac{2}{\lambda_{\rm min}(\bfB_*)} \left( \sqrt{\frac{4C_2^2 \ln((p+1)/\delta)}{M}} + \frac{4C_2 \ln((p+1)/\delta)}{3M} \right) \\
    & + \frac{2}{\lambda_{\rm min}^2(\bfB_*)} \left( \sqrt{\frac{4C_1^2 \ln(2p/\delta)}{M}} + \frac{4C_1 \ln(2p/\delta)}{3M} \right) C_2 \\
    & = \mathcal{O} (\frac{1}{\sqrt{M}}).
\end{align}

\end{document}